\def\Ass{\mathop{\rm Ass}\nolimits}

\def\End{\mathop{\rm End}\nolimits}

\def\height{\mathop{\rm height}\nolimits}
\def\Hom{\mathop{\rm Hom}\nolimits}

\def\dirlim{\mathop{\rm dirlim}\nolimits}
\def\invlim{\mathop{\rm invlim}\nolimits}

\def\InjH{\mathop{\rm E}\nolimits}
\def\LCMo{\mathop{\rm H}\nolimits}

\def\Spec{\mathop{\rm Spec}\nolimits}

\def\char{\mathop{\rm char}\nolimits}
\def\Naturalsign{{\rm l\kern-.23em N}}
\baselineskip=15pt \font \normal=cmr10 scaled \magstep0 \font
\mittel=cmr10 scaled \magstep3 
\input amssym.def
\input amssym.tex
{\parindent=0pt \mittel Finiteness properties of duals of local
cohomology modules \normal
\bigskip Michael Hellus
\footnote{(+)}{The author was partially supported by the DFG}
\par
E-mail: michael.hellus@math.uni-leipzig.de
\bigskip
Abstract:
\par
{\tt We investigate Matlis duals of local cohomology modules and
prove that, in general, their \par zeroth Bass number with respect
to the zero ideal is not finite. We also prove that,\par somewhat
surprisingly, if we apply local cohomology again (i. e. to the
Matlis dual of \par the local cohomology module), we get (under
certain hypotheses) either zero or $\InjH $, an $R$- \par injective
hull of the residue field of the local ring $R$.}
\bigskip
{\bf 0 Introduction}} \smallskip Let $I$ be an ideal of a noetherian
local ring $(R,\goth m)$; $\LCMo ^i_I(R)$ is the $i$-th local
cohomology module of $R$ supported in $I$ (see [Gr], [BS] for
general facts on these modules) and $\InjH :=\InjH _R(R/\goth m)$ a
fixed $R$-injective hull of $R/\goth m$. By $D$ we denote the Matlis
dual functor from $(R-mod)$ to $(R-mod)$, i. e.
$$D(M)=\Hom _R(M,\InjH )$$
for every $R$-module $M$.
\par
Recently there was some work on the modules $D(\LCMo ^i_I(R))$, see
e. g. [HS1], [HS2] and [H1] -- [H5]; one major motivation for the
study of these modules is remark 1.1 below (originally [H5, Corollay
1.1.4] or [H1, section 0]) which connects regular sequences on
$D(\LCMo ^i_I(R))$ in a clear way to the notion of a set-theoretic
complete intersection ideal (i. e. an ideal which is, up to radical,
generated by a regular sequence). Besides remark 1.1 various
applications of results on $D(\LCMo ^i_I(R))$ are collected in [H5,
section 6].
\par
In some cases it is known (see e. g. [H5, theorem 3.1.3]) that the
zero ideal of $R$ (assume $R$ is a domain) is associated to the
$R$-module $D(\LCMo ^i_I(R))$; there is a conjecture (*) on the set
of associated prime ideals of $D(\LCMo ^i_{\underline x}(R))$, where
$\underline x=x_1,\dots ,x_i$ is a sequence in $R$: It is conjecture
(*) from [H1] and [H5]:
$$\Ass _R(D(\LCMo ^i_{(x_1,\dots ,x_i)R}(R)))=\{ \goth p\in
\Spec (R)\vert \LCMo ^i_{(x_1,\dots ,x_i)R}(R/\goth p)\neq 0\} \ \
.\leqno{(*)} $$ Clearly, if conjecture (*) holds, we have
$$\{ 0\} \in \Ass _R(D(\LCMo ^i_{\underline x}(R)))$$
provided $\LCMo ^i_{\underline x}(R)\neq 0$ holds. More details on
conjecture (*) are contained in [H1] and [H5]. In any case, it is
natural to ask for the associated Bass number, i. e. for the
$Q(R)$-vector space dimension of
$$D(\LCMo ^i_I(R))\otimes _RQ(R)\ \ ,$$
where $Q(R)$ is the quotient field of $R$. Theorem 2.3, which is the
first main result of this work, shows that this number is not
finite, in general. In contrast to this, the results in section 3
show that $D(\LCMo ^i_I(R))$ is "small" in the following sense:
$\LCMo ^i_I(D(\LCMo ^i_I(R)))$ is either $\InjH $ or zero; more
precisely, it is $\InjH $ if $I$ is a set-theoretic complete
intersection ideal (theorem 3.1), while it is zero or $\InjH $ in a
more general situation (theorem 3.2); note that, for theorem 3.2, we
use the so-called $D$-module structure on $\LCMo ^i_I(D(\LCMo
^i_I(R)))$, in particular we use [Ly, theorem 2.4] in connection
with zero-dimensional $D$-modules and thus we have to assume
equicharacteristic zero in the statement of theorem 3.2.
\par
The author thanks Gennady Lyubeznik for helpful discussions, in
particular on $D$-modules. {\parindent=0pt
\bigskip
{\bf 1 Prerequisites} \smallskip This section collects some
statements on Matlis duals of local cohomology modules that will be
needed in the next two sections. For instance, the following is
known about regular sequences on such Matlis duals:
\bigskip
{\bf 1.1 Remark ([H5, Corollary 1.1.4])}
\par
Let $(R,\goth m)$ be a noetherian local ring, $I$ a proper ideal of
$R$, $h\in \Naturalsign $ and $\underline f=f_1,\dots ,f_h\in I$ an
$R$-regular sequence. The following statements are equivalent:
\par
(i) $\sqrt {\underline fR}=\sqrt I$ (in this case $I$ is a
set-theoretic complete intersection).
\par
(ii) $\LCMo ^l_I(R)=0$ for every $l>h$ and the sequence $\underline
f$ is quasi-regular on $D(\LCMo ^h_I(R))$.
\par
(iii) $\LCMo ^l_I(R)=0$ for every $l>h$ and the sequence $\underline
f$ is regular on $D(\LCMo ^h_I(R))$.
\par
(The case $h=0$ means
$$\eqalign {\sqrt I=\sqrt 0&\iff \LCMo ^l_I(R)=0 \hbox{ for every }l>0\cr &\iff \LCMo ^l_I(R)=0 \hbox{ for every }l>0 \hbox { and } \Gamma _I(R)\neq 0\ \ ).\cr }$$
{\bf 1.2 Remark}
\par
We want to calculate local cohomology modules and their Matlis duals
in the following situation: $k$ a field, $R=k[[X_1,\dots ,X_n]]$ a
power series algebra over $k$ in $n\in \Naturalsign $ variables,
$i\in \{ 0,\dots ,n\} $ and $I$ the ideal $(X_1,\dots ,X_i)R$ of
$R$, we want to calculate $\LCMo ^i_I(R)$ and $D(\LCMo ^i_I(R))$:
\smallskip
By considering the \v{C}ech-complex of $R$ with respect to
$X_1,\dots ,X_i$ it is not difficult to see that there is a
canonical equality
$$\LCMo ^i_I(R)=k[[X_{i+1},\dots ,X_n]][X_1^{-1},\dots ,X_i^{-1}]$$
(by using this notation we mean, of course, the free
$k[[X_{i+1},\dots ,X_n]]$-module on the set of "inverse monomials"
$\{ X_1^{\nu _1}\cdot \dots \cdot X_i^{\nu _i}\vert \nu _1,\dots
,\nu _i\leq 0\} $; this module has an obvious $R$-module structure).
It is well-known and easy to see that one can realize $\InjH :=\InjH
_R(k)$ as $k[X_1^{-1},\dots ,X_n^{-1}]$ (like before, this is the
$k$-vector space with basis $\{ X_1^{\nu _1}\cdot \dots \cdot
X_n^{\nu _n}\vert \nu _1,\dots ,\nu _n\leq 0\} $ having a natural
$R$-module structure and containing $k=k\cdot 1=k\cdot X_1^0\cdot
\dots \cdot X_n^0$).\smallskip It is straight-forward to see that
every element of
$$k[X_{i+1}^{-1},\dots ,X_n^{-1}][[X_1,\dots
,X_i]]$$ (i. e. every formal power series in $X_1,\dots ,X_i$ and
coefficients in $k[X_{i+1}^{-1},\dots ,X_n^{-1}]$) determines
(essentially by multiplication) an $R$-linear map
$$k[[X_{i+1},\dots ,X_n]][X_1^{-1},\dots ,X_i^{-1}]\to k[X_1^{-1},\dots
,X_n^{-1}]\ \ ,$$ i. e. an element of $D(\LCMo ^i_I(R))$.
Furthermore, the element of $k[X_{i+1}^{-1},\dots
,X_n^{-1}][[X_1,\dots ,X_i]]$ is uniquely determined by its
associated map $\LCMo ^i_I(R)\to \InjH $; and finally, a tedious but
easy calculation shows that every $R$-linear map $\LCMo ^i_I(R)\to
\InjH $ arises this way, i. e. there is an equality
$$D(\LCMo ^i_I(R))=k[X_{i+1}^{-1},\dots ,X_n^{-1}][[X_1,\dots
,X_i]]\ \ .$$ {\bf 1.3 Remark (see [H5, subsection 7.2])}
\par
Let $k$ be a field and $R=k[[X_1,\dots X_n]]$ a power series ring
over $k$ in $n$ variables. Let
$$D(R,k)\subseteq \End _k(R)$$
be the (non-commutative) subring defined by the multiplication maps
by $r\in R$ (for all $r\in R$) and by all $k$-linear derivation maps
from $R$ to $R$. $D(R,k)$ is the so-called ring of $k$-linear
differential operators on $R$. [Bj] contains material on the ring
$D(R,k)$ and on similar rings; $D(R,k)$-modules in relation with
local cohomology modules have been studied in [Ly]. For $i=1,\dots
,n$ let $\partial _i$ denote the partial derivation map from $R$ to
$R$ with respect to $X_i$. Then, as an $R$-module, one has
$$D(R,k)=\oplus _{i_1,\dots
,i_n\in \Naturalsign } R\cdot \partial _1^{i_1}\dots \partial
_n^{i_n}\ \ .\leqno{(1)}$$ Now, let $I\subseteq R$ be an ideal and
$i\in \Naturalsign $. We will demonstrate that there is a canonical
left-$D(R,k)$-module structure on $D(\LCMo ^i_I(R))$. To do so, by
identity (1), it is sufficient to determine the action of an
arbitrary $k$-linear derivation $\delta :R\to R$ on $D(\LCMo
^i_I(R))$, to extend it to an action of $D(R,k)$ on $D(\LCMo
^i_I(R))$ and to show that this action is well-defined and satisfies
all axioms of a left-$D(R,k)$-module. The derivation $\delta $
induces a $k$-linear map
$$R/I^v\to R/I^{v-1}\ \ (v\geq 1)$$
and, in a canonical way, a map of complexes from the \v{C}ech
complex of $R/I^v$ with respect to $X_1,\dots ,X_n$ to the \v{C}ech
complex of $R/I^{v-1}$ with respect to $X_1,\dots ,X_n$ ($v\geq 1$).
By taking cohomology, we get a map
$$\LCMo ^{n-i}_\goth m(R/I^v)\to \LCMo ^{n-i}_\goth m(R/I^{v-1})\ \
(v\geq 1)\ ,$$ where $\goth m$ stands for the maximal ideal of $R$.
These maps induce a map
$$\invlim _{v\in \Naturalsign }(\LCMo ^{n-i}_\goth m(R/I^v))\to \invlim _{v\in \Naturalsign }(\LCMo ^{n-i}_\goth
m(R/I^v))$$ (note that the maps of the above inverse system are
induced by the canonical epimorphisms $R/I^v\to R/I^{v-1}$). But, by
local duality and $\LCMo ^i_I(R)=\dirlim _{v\in \Naturalsign }(Ext
^i_R(R/I^v,R))$, one has
$$\invlim _{v\in \Naturalsign }(\LCMo ^{n-i}_\goth m(R/I^v))=D(\dirlim
_{v\in \Naturalsign }(Ext ^i_R(R/I^v,R)))=D(\LCMo ^i_I(R))\ \ .$$
Now, having determined the action of the element $\delta $ on
$D(\LCMo ^i_I(R))$, by (1) it is clear how to extend this to an
action of $D(R,k)$ on $D(\LCMo ^i_I(R))$ such that $D(\LCMo
^i_I(R))$ becomes a left-$D(R,k)$-module (note that, for every
$k$-linear derivation $\delta :R\to R$ and every $r\in R$, we have
$\delta (r\cdot d)=\delta (r)\cdot d+r\cdot \delta (d)$, i. e. the
action of $D(R,k)$ on $D(\LCMo ^i_I(R)$ makes it a
left-$D(R,k)$-module). It is known (e. g. from various results in
[H5, sections 2 -- 4]) that, in general,
$$D(\LCMo ^i_I(R))$$
has infinitely many associated primes. On the other hand, one knows
from [Ly, theorem 2.4 (c)] (at least if $\char (k)=0$), that every
finitely generated left-$D(R,k)$-module has only finitely many
associated prime ideals (as $R$-module, of course). This shows that,
in general, $D(\LCMo ^i_I(R))$ is an example of a non-finitely
generated left-$D(R,k)$-module. In particular, $D(\LCMo ^i_I(R))$ is
not holonomic in general (see [Bj] for the notion of holonomic
modules).
\bigskip
{\bf 2 The zeroth Bass number of $D(\LCMo ^i_I(R))$ (w. r. t. the
zero ideal) is not finite in general} \smallskip Let $(R,\goth m)$
be a noetherian local domain, $i\geq 1$ and $x_1,\dots ,x_i\in R$.
Then one has
$$\{ 0\} \in \Ass _R(D(\LCMo ^i_{(x_1,\dots ,x_i)R}(R)))$$
in some situations (see e. g. [H5, theorem 3.1.3]); actually, if
conjecture (*) (from [H1] and [H5]) holds, this is true provided
$\LCMo ^i_{(x_1,\dots ,x_i)R}(R)\neq 0$ holds. It is natural to ask
for the associated Bass number of $D(\LCMo ^i_{(x_1,\dots
,x_i)R}(R))$, i. e. the $Q(R)$-vector space dimension of
$$D(\LCMo ^i_{(x_1,\dots ,x_i)R}(R))\otimes _RQ(R)\ \ ,$$
where $Q(R)$ stands for the quotient field of $R$. As we will see
below, this number is not finite in general; more precisely, we
consider the following case: Let $k$ be a field, $R=k[[X_1,\dots
,X_n]]$ a power series algebra over $k$ in $n\geq 2$ variables,
$1\leq i<n$ and $I$ the ideal $(X_1,\dots ,X_i)R$ of $R$; in this
situation
$$\dim _{Q(R)}(D(\LCMo ^i_I(R))\otimes _RQ(R))=\infty $$
holds; this is the statement of theorem 2.3. \bigskip {\bf 2.1
Remark}
\par
Note that in section 1.2 we introduced some notation on polynomials
in "inverse variables" and we explained and proved the following
formulas:
$$\LCMo ^i_I(R)=k[[X_{i+1},\dots ,X_n]][X_1^{-1},\dots ,X_i^{-1}]\ \ ,$$
$$\InjH _R(k)=k[X_1^{-1},\dots ,X_n^{-1}]$$
and
$$D(\LCMo ^i_I(R))=k[X_{i+1}^{-1},\dots ,X_n^{-1}][[X_1,\dots
,X_i]]\ \ .$$ Also note that the latter module is different from
(and larger) than the module
$$k[[X_1,\dots ,X_i]][X_{i+1}^{-1},\dots ,X_n^{-1}]\ \ .$$
\bigskip
{\bf 2.2 Remark} \par The following proof of theorem 2.3, which says
that a certain Bass number is infinite, is technical; its basic idea
is the following one: Let $k$ be a field, $R=k[[X,Y]]$ a power
series algebra over $k$ in two variables; then we have
$$\LCMo ^1_{XR}(R)=k[[Y]][X^{-1}]$$
and $$D:=D(\LCMo ^1_{XR}(R))=k[Y^{-1}][[X]]\ \ .$$ Set
$$\eqalign {d_2&:=\sum _{l\in \Naturalsign }Y^{-l^2}X^l\cr
&=1+Y^{-1}X+Y^{-4}X^2+Y^{-9}X^3+\dots \in D\cr }$$ and let $r\in
R\setminus \{ 0\} $ be arbitrary. Because of $r\neq 0$ we can write
$$r=X^{a+1}\cdot h+X^a\cdot g$$
with some $h\in R, g\in k[[Y]]\setminus \{ 0\} $. Then, at least for
$l>>0$, the coefficient of $r\cdot d_2$ in front of $X^l$ is
$$h^*\cdot Y^{-(l-a-1)^2}+g\cdot Y^{-(l-a)^2}$$
for some $h^*\in k[[Y]]$. Now, if we write
$$g=c_bY^b+c_{b+1}Y^{b+1}+\dots $$
for some $b\in \Naturalsign , c_b\neq 0$ and observe the fact
$$-(l-a)^2+b<-(l-a-1)^2 \ \ (l>>0)\ \ ,$$
it follows that the term
$$c_b\cdot Y^{-(l-a)^2+b}$$
(coming from $h^*\cdot Y^{-(l-a-1)^2}+g\cdot Y^{-(l-a)^2}$) cannot
be canceled out by any other term. In fact, for $l>>0$, the lowest
non-vanishing $Y$-exponent of the coefficient in front of $X^l$, is
$-(l-a)^2+b$. The crucial point is that the sequences $-(l-a)^2+b$
and $-l^2$ agree up to the two shifts given by $a$ and $b$. This
means that some information about $d_2$ is stored in $rd_2$.
\bigskip
{\bf 2.3 Theorem}
\par
Let $k$ be a field, $R=k[[X_1,\dots ,X_n]]$ a power series algebra
over $k$ in $n\geq 2$ variables, $1\leq i<n$ and $I$ the ideal
$(X_1,\dots ,X_i)R$ of $R$. Then
$$\dim _{Q(R)}(D(\LCMo ^i_I(R))\otimes _RQ(R))=\infty $$
holds.
\par
Proof:
\par
As the proof is technical we will first show the case $n=2,i=1$; in
the remark after this proof we will explain how one can reduce the
general to this special case. Set $X=X_1,Y=X_2$ and
$$D:=D(\LCMo ^2_I(R))=k[Y^{-1}][[X]]\ \ .$$
For every $n\in \Naturalsign \setminus \{ 0\} $, set
$$d_n:=\sum _{l\in \Naturalsign }Y^{-l^n}\cdot X^l\in D\ \ .$$
It is sufficient to show the following statement: The elements
$(d_n\otimes 1)_{n\in \Naturalsign \setminus \{ 0\} }$ in $D\otimes
_RQ(R)$ are $Q(R)$-linear independent: \smallskip We define an
equivalence relation on ${\bf Z}^\Naturalsign $ (the set of all maps
from $\Naturalsign $ to $\bf Z$, i. e. infinite sequences of
integers) by saying that $(a_n),(b_n)\in {\bf Z}^\Naturalsign $ are
equivalent (short form: $(a_n)\sim (b_n)$) iff there exist $N,M\in
\Naturalsign $ and $p\in {\bf Z}$ such that
$$a_{N+1}=b_{M+1}+p,a_{N+2}=b_{M+2}+p,\dots $$
hold. It is easy to see that $\sim $ is an equivalence relation on
${\bf Z}^\Naturalsign $. For every $d\in D$, we define $\delta
(d)\in {\bf Z}^\Naturalsign $ in the following way: Let $f_l\in
k[Y^{-1}]$ be the coefficient of $d$ in front of $X^l$; we set
$$(\delta (d))(l):=0$$
if $f_l=0$ and
$$(\delta (d))(l):=s$$
if $s$ is the smallest $Y$-exponent of $f_l$, i. e.
$$f_l=c_sY^s+c_{s+1}Y^{s+1}+\dots +c_0\cdot 1$$
for some $c_s\neq 0$.
\bigskip
Now suppose that $r_1,\dots ,r_{n_0}\in R$ are given such that
$r_{n_0}\neq 0$. We claim that
$$\delta (r_1d_1+\cdot +r_{n_0}d_{n_0})\sim \delta (d_{n_0})$$
holds. Note that if we prove this statement we are done, essentially
because then $r_1d_1+\dots +r_{n_0}d_{n_0}$ can not be zero.\bigskip
It is obvious that one has $\delta (d+d^\prime )\sim \delta
(d_{N_2})$ for given $d,d^\prime \in D$ such that
$$\delta (d)\sim \delta (d_{N_1}), \delta (d^\prime )\sim \delta
(d_{N_2}), N_2>N_1\ \ .$$ For this reason it is even sufficient to
prove the following statement: For a fixed $n\in \Naturalsign
\setminus \{ 0\} $ and for any $r\in R\setminus \{ 0\} $ one has
$$\delta (rd_n)\sim \delta (d_n)\ \ .$$
We can write
$$r=X^{a+1}\cdot h+X^a\cdot g$$
with $a\in \Naturalsign , h\in k[[X,Y]]$ and $g\in k[[Y]]\setminus
\{ 0\} $. We get
$$\delta (r\cdot d_n)\sim \delta (\sum _{l\geq a+1}
(hY^{-(l-a-1)^n}+gY^{-(l-a)^n})X^l)$$ and we write
$$g=c_bY^b+c_{b+1}Y^{b+1}+\dots $$
with $c_b\in k^*$. Now, because of
$$-(l-a)^n+b<-(l-a-1)^n\ \ (l>>0)$$
it is clear that, for $l>>0$, the smallest $Y$-exponent in front of
$X^l$ (of the power series $r\cdot d_n$) is $-(l-a)^n+b$. Therefore,
one has
$$\delta (r\cdot d_n)\sim (-l^n)\sim \delta (d_n)$$
and we are done.
\bigskip
{\bf 2.4 Remark}
\par
A proof of the general case of theorem 2.3 can be obtained e. g. in
the following way: First, we use [H5, theorem 3.1.2] repeatedly to
get a surjection
$$\LCMo ^i_{(X_1,\dots ,X_i)R}(R)\to \LCMo ^{n-1}_{(X_1,\dots
,X_{n-1})R}(R)$$ and hence an injection
$$D(\LCMo ^{n-1}_{(X_1,\dots
,X_{n-1})R}(R))\to D(\LCMo ^i_{(X_1,\dots ,X_i)R}(R))\ \ ,$$ which
allows us to reduce to the case $i=n-1$; then it is possible to
adapt our proof of theorem 2.3 with some minor changes: Instead of
working with maps $\Naturalsign \to {\bf Z}$, one works with maps
$$\Naturalsign ^{n-1}\to {\bf Z}$$
and also with multi-indices instead of indices.
\bigskip
{\bf 3 On the module $\LCMo ^i_I(D(\LCMo ^i_I(R)))$} \smallskip In
this section we assume $\LCMo ^l_I(R)=0$ for every $l>\height (I)$
and investigate the local cohomology module
$$\LCMo ^{\height (I)}_I(D(\LCMo ^{\height (I)}_I(R)))\ \ .$$
Our results say (essentially) that this module is $\InjH _R(R/\goth
m)$ if $I$ is a set-theoretic complete intersection (theorem 3.1)
and that it is either $\InjH _R(R/\goth m)$ or zero in general
(theorem 3.2):
\bigskip
{\bf 3.1 Theorem}
\par
Let $(R,\goth m)$ be a noetherian local complete Cohen-Macaulay ring
with coefficient field $k$ and $x_1,\dots ,x_i\in R$ ($i\geq 1$) a
regular sequence in $R$. Set $I:=(x_1,\dots ,x_i)R$ ($I$ is a
set-theoretic complete intersection ideal of $R$). Then one has
$$\LCMo ^i_I(D(\LCMo ^i_I(R)))=\InjH _R(k)\ \ .$$
Proof:
\par
First we show a special case: Assume that $R=k[[X_1,\dots ,X_n]]$ is
a formal power series algebra over $k$ in $n$ variables and
$x_1=X_1,\dots ,x_i=X_i$. Then, as we have seen in remark 1.3, we
can write
$$\LCMo ^i_I(R)=k[[X_{i+1},\dots ,X_n]][X_1^{-1},\dots ,X_i^{-1}]$$
and
$$D(\LCMo ^i_I(R))=k[X_{i+1}^{-1},\dots ,X_n^{-1}][[X_1,\dots
,X_i]]\ \ .$$
As the functor $\LCMo ^i_I$ is right-exact, we have
$$\eqalign {\LCMo ^i_I(D(\LCMo ^i_I(R)))&=\LCMo ^i_I(R)\otimes _RD(\LCMo
^i_I(R))\cr &=k[[X_{i+1},\dots ,X_n]][X_1^{-1},\dots
,X_i^{-1}]\otimes _Rk[X_{i+1}^{-1},\dots X_n^{-1}][[X_1,\dots
,X_i]]\cr &\buildrel (*)\over =k[X_1^{-1},\dots ,X_n^{-1}]\cr
&=\InjH _R(k)\cr }$$ Proof of equality (*): The map
$$k[X_{i+1},\dots ,X_n][X_1^{-1},\dots
,X_i^{-1}]\otimes k[X_{i+1}^{-1},\dots X_n^{-1}][X_1,\dots ,X_i]\to
k[X_1^{-1},\dots ,X_n^{-1}]$$
$$X_{i+1}^{r_{i+1}}\cdot \dots \cdot X_n^{r_n}\cdot X_1^{-s_1}\cdot
\dots X_i^{-s_i}\otimes X_{i+1}^{-t_{i+1}}\cdot \dots \cdot
X_n^{-t_n}\cdot X_1^{u_1}\cdot \dots \cdot X_i^{u_i}\mapsto $$
$$\mapsto X_{i+1}^{r_{i+1}-t_{i+1}}\cdot X_n^{r_n-t_n}\cdot
X_1^{u_1-s_1}\cdot \dots \cdot X_i^{u_i-s_i}\hbox { if
}r_{i+1}-t_{i+1},\dots ,r_n-t_n,u_1-s_1,\dots ,u_i-s_i\leq 0$$ and
to zero otherwise, induces an $R$-linear map
$$k[[X_{i+1},\dots ,X_n]][X_1^{-1},\dots
,X_i^{-1}]\otimes _Rk[X_{i+1}^{-1},\dots X_n^{-1}][[X_1,\dots
,X_i]]\to k[X_1^{-1},\dots ,X_n^{-1}]\ \  ,\leqno {(1)}$$ which is
surjective and maps the $k$-vector space generating system
$$\{ X_1^{-s_1}\cdot \dots \cdot X_i^{-s_i}\otimes
X_{i+1}^{-t_{i+1}}\cdot \dots \cdot X_n^{-t_n}\vert s_1,\dots
,s_i,t_{i+1},\dots ,t_n\geq 0\} $$ of the vector space on the left
side of (1) to the $k$-basis
$$\{ X_1^{-s_1}\cdot \dots \cdot X_i^{-s_i}\cdot
X_{i+1}^{-t_{i+1}}\cdot \dots \cdot X_n^{-t_n}\vert s_1,\dots
,s_i,t_{i+1},\dots ,t_n\geq 0\} $$ of the vector space on the right
side of (1), and therefore provides us with the desired isomorphism
in our special case.
\par
We come to the general case: Choose $x_{i+1},\dots ,x_n\in R$ such
that $\sqrt {(x_1,\dots ,x_n)R}=\goth m$ ($x_1,\dots x_n$ is a s. o.
p. of $R$). Define
$$R_0:=k[[x_1,\dots ,x_n]]\subseteq R$$
$R_0$ is regular of dimension $n$ and $R$ is a finite-rank free
$R_0$-module. Define $I_0:=(x_1,\dots ,x_i)R_0$. We have
$$\LCMo ^i_I(D(\LCMo ^i_I(R)))=\LCMo ^i_I(R)\otimes _RD(\LCMo
^i_I(R))$$ and
$$\LCMo ^i_I(R)=\LCMo ^i_{I_0}(R_0)\otimes _{R_0}R$$
and
$$\eqalign {D(\LCMo ^i_I(R)))&=\Hom _R(\LCMo ^i_{I_0}(R_0)\otimes
_{R_0}R,\InjH _R(k))\cr &=\Hom _{R_0}(\LCMo ^i_{I_0}(R_0),\InjH
_R(k))\cr &\buildrel (2)\over =\Hom _{R_0}(\LCMo ^i_{I_0}(R_0),\Hom
_{R_0}(R,\InjH _{R_0}(k)))\cr &=\Hom _{R_0}(R,D_{R_0}(\LCMo
^i_{I_0}(R_0)))\cr }$$ For (2) we use the fact
$$\InjH _R(k)=\Hom _{R_0}(R,\InjH _{R_0}(k))$$
We get
$$\eqalign {\LCMo ^i_I(D(\LCMo ^i_I(R)))&=\LCMo ^i_{I_0}(R_0)\otimes _{R_0}\Hom _{R_0}(R,D_{R_0}(\LCMo
^i_{I_0}(R_0)))\cr &\buildrel (3)\over =\Hom _{R_0}(R,\LCMo
^i_{I_0}(R_0)\otimes _{R_0}D_{R_0}(\LCMo ^i_{I_0}(R_0)))\cr &=\Hom
_{R_0}(R,\InjH _{R_0}(k))\cr &\buildrel (2)\over =\InjH _R(k)\cr }$$
For (3) we use the fact that $R$ is a finite-rank free $R_0$-module.
\vfil \eject {\bf 3.2 Theorem}
\bigskip
Let $R$ be a noetherian local complete regular ring of
equicharacteristic zero, $I\subseteq R$ an ideal of height $i\geq
1$, $x_1,\dots ,x_i\in I$ an $R$-regular sequence and assume that
$$\LCMo ^l_I(R)=0\hbox { for every }l>i\ \ .$$
Then $\LCMo ^i_I(D(\LCMo ^i_I(R)))$ is either $\InjH _R(k)$ or zero.
\par
Proof:
\par
We set
$$D:=D(\LCMo ^i_{(x_1,\dots ,x_i)R}(R))\ \ $$
By remark 1.1, we know that $x_1,\dots ,x_i$ is a $D$-regular
sequence and, therefore, we have
$$\LCMo ^0_{(x_1,\dots ,x_i)R}(D)=\dots =\LCMo ^{h-1}_{(x_1,\dots
,x_i)R}(D)=0\ \ .$$ Because of this, an easy spectral sequence
argument (applied to the composed functor $\Gamma _I\circ \Gamma
_{(x_1,\dots ,x_i)R}$ and to the $R$-module $D$) shows that
$$\LCMo ^i_I(D)=\Gamma _I(\LCMo ^i_{(x_1,\dots ,x_i)R}(D))\subseteq
\LCMo ^i_{(x_1,\dots ,x_i)R}(D)=\InjH _R(k)\ \ .$$ The last equality
is theorem 3.1. But, from remark 1.3 and from [Ly, Example 2.1
(iv)], it is clear that
$$\LCMo ^i_I(R), D(\LCMo ^i_I(R))\hbox{ and } \LCMo ^i_I(D(\LCMo ^i_I(R)))$$
all have a $D(R,k)$-module structure and so, from [Ly, theorem 2.4
(b)], we deduce that $\LCMo ^i_I(D)$ is either $\InjH _R(k)$ or
zero. Furthermore, the natural injection
$$\LCMo ^i_I(R)\subseteq \LCMo ^i_{(x_1,\dots ,x_i)R}(R)$$
induces a surjection
$$D\to D(\LCMo ^i_I(R))$$
and hence, as $\LCMo ^i_I$ is right-exact, a surjection
$$\LCMo ^i_I(D)\to \LCMo ^i_I(D(\LCMo ^i_I(R)))\ \ .$$
But again, the last module has a $D(R,k)$-module structure, and
thus, from [Ly, theorem 2.4 (b)] and from what we know already, we
conclude the statement. \bigskip \bigskip {\bf References} \normal
\bigskip
\parindent=1.2cm
\def\litem{\par\noindent \hangindent=\parindent\ltextindent}
\def\ltextindent#1{\hbox to \hangindent{#1\hss}\ignorespaces}
\litem{[BH]} Bruns, W. and Herzog, J. Cohen-Macaulay Rings, {\it
Cambridge University Press}, (1993).
\medskip
\litem{[Bj]} Bjork, J.-E. Rings of Differential Operators, {\it
Amsterdam North-Holland}, (1979)
\medskip
\litem{[BS]} Brodmann, M. P. and Sharp, R. J. Local Cohomology, {\it
Cambridge studies in advanced mathematics} {\bf 60}, (1998).
\medskip
\litem{[Gr]} Grothendieck, A. Local Cohomology, {\it Lecture Notes
in Mathematics, Springer Verlag}, (1967).
\medskip
\litem{[H1]} Hellus, M. On the associated primes of Matlis duals of
top local cohomology modules, {\it Communications in Algebra} {\bf
33}, (2001), no. 11, 3997--4009.
\medskip
\litem{[H2]} Hellus, M. Matlis duals of top local cohomology modules
and the arithmetic rank of an ideal, to appear in {\it
Communications in Algebra}.
\medskip
\litem{[H3]} Hellus, M. Attached primes and Matlis duals of local
cohomology modules, submitted to {\it Archiv der Mathematik}
\medskip
\litem{[H4]} Hellus, M. Local Homology, Cohen-Macaulayness and
Cohen-Macaulayfications, to appear in {\it Algebra Colloquium}.
\medskip
\litem{[H5]} Hellus, M. Local Cohomology and Matlis duality, {\it
Habilitationsschrift}, Leipzig (2006).
\medskip
\litem{[HS1]} Hellus, M. and St\"uckrad, J. Matlis duals of top
Local Cohomology Modules, submitted to {\it Proceedings of the
American Mathematical Society}.
\medskip
\litem{[HS2]} Hellus, M. and St\"uckrad, J. Generalization of an
example of Hartshorne concerning local cohomology, preprint.
\medskip
\litem{[Ly]} Lyubeznik, G. Finiteness properties of local cohomology
modules (an application of $D$-modules to Commutative Algebra), {\it
Invent. Math.} {\bf 113}, (1993), 41--55.
\medskip
\litem{[Ms]} Matlis, E. Injective modules over Noetherian rings,
{\it Pacific J. Math.} {\bf 8}, (1958) 511--528. }\medskip
\end